\documentclass{article}
\usepackage{amsmath}

\begin{document}

\begin{center}
\textbf{Phase flows and vectorial lagrangians in $J^3(\pi)$}

\bigskip

V.N.Dumachev \\
Voronezh institute of the MVD of Russia\\
e-mail: dumv@comch.ru
\end{center}

\textbf{Abstract.} On the basis of Liouville theorem the
generalization of the Nambu mechanics is considered.  For
three-dimensional phase space the concept of vector hamiltonian and
vector lagrangian is entered.

\bigskip
\textbf{1.} Standard phenomenological approach to the analysis
dynamic system is the construction for it the functional of actions
$S=\int L\;dt$. We represent this functional as submanifolds in jet
bundles $J^n(\pi)$: $E \to M$
\[
 F(t,x_0,x_1,...,x_n)=0,
\]
where $t \in M \subset R$, $u = x_{0} \in U \subset R$, $x_i \in
J^i(\pi) \subset R^n$, $E=M \times U$.

The Euler-Lagrange equation
\begin{equation}
\label{eq1} \sum\limits_{k=0}^n(-)^k\frac{d^k}{dt^k}\frac{\partial
L}{\partial x_k}=0
\end{equation}
describes a lines (jet) in $J^{2n}(\pi)$. Embedding $J^n(\pi)
\subset J^1( J^1( J^1(... J^1(\pi)))))$ allows us to rewrite
differential equation $n$-st order as system of the $n$ equations
1-st order
\begin{equation}
\label{eq2} \overset{\cdot}{\textbf{x}}=\textbf{Ax}.
\end{equation}
According to the Noether theorem, symmetry of functional $S$ with
respect to generator $X=\partial/\partial t$ give us the
conservation law $I$, and hamiltonian form for our dynamics:
\begin{equation}
\label{eq3} \overset{\cdot}{\textbf{x}}=\{H(I),\textbf{x}\}.
\end{equation}

\bigskip
\textbf{2.} Another approach for receiving of the Euler-Lagrange
equation (\ref{eq1}) for every hamiltonians set was described by
P.A.Griffiths \cite{Griffiths}. He find such 1-form
\[
\psi=Ldt+\lambda^i \theta_i, \qquad i=0..n-1,
\]
which does not vary at pullback along vector fields
$(\partial/\partial \theta_i,\partial/\partial d\theta_{n-1})$. Here
\[
\theta_i=dx_i-x_{i+1}\; dt
\]
is the contact distribution, $\lambda^i$ is the Lagrange
multipliers. Bounding of the form $\Psi=d\psi$ on the field
$(\partial/\partial \theta_i,\partial/\partial d\theta_{n-1})$ gives
the set
\[
\left\{
\begin{array}{l}
\dfrac{\partial L}{\partial
x_i}=\dfrac{d\lambda^i}{dt}+\lambda^{i-1},
\qquad i=0..n-1,\\
\\
\dfrac{\partial L}{\partial x_n}=\lambda^{n-1},
\end{array}\right.
\]
which is equivalent to the equation (\ref{eq1}).

Hamiltonian formulation of this theory assume that Lagrange's
multipliers $\lambda^i$ be a dynamic variables $H=H(x_i,\lambda^i)$:
\[
\psi=(L-\lambda^i x_{i+1})\wedge dt+\lambda^i dx_i=-H dt+\lambda^i
dx_i.
\]
Then bounding of the form $\Psi=d\psi$ on the field
$(\partial/\partial x_i,\partial/\partial \lambda^i)$ gives
\[
\frac{\partial H}{\partial x_i}=-\frac{\partial \lambda^i}{\partial
t}, \qquad \frac{\partial H}{\partial \lambda^i}=\frac{\partial
x_i}{\partial t}.
\]

\bigskip
\textbf{3.} Our target is the generalization of the above scheme on
a case odd jets. To clear idea of a method we shall receive the
Euler-Lagrange equation for $L \in J^3(\pi)$.

\textbf{Theorem 1.} For $L \in J^3(\pi)$ the Euler-Lagrange equation
has the form
\begin{equation}
\label{eq4} \frac{1}{2}\frac{d}{dt}\left( L_{\overset{\cdot
}{x}_{k}}^{i}-L_{\overset{\cdot }{x}_{i}}^k\right) =
L^k_{x_i}-L^i_{x_k}.
\end{equation}

\textbf{Proof.} Let $\psi$ be the Griffiths 2-form:
\[
\psi=L^idx_i\wedge dt+\lambda^i \Theta_i,
\]
where $\Theta = \theta \wedge \theta$. Exterior differencial this
form is
\begin{eqnarray*}
d\psi &=&dL^{i}\wedge \omega _{i}=\left( \text{rot}\,L\right) ^{k}\Theta _{k}\wedge dt+d\lambda ^{i}\wedge \Theta _{i}+\lambda ^{i}\wedge d\Theta _{i}+ \\
&+&L_{\overset{\cdot }{x}_{3}}^{2}d\overset{\cdot }{x}_{3}\wedge
dx_{2}\wedge dt+L_{\overset{\cdot }{x}_{2}}^{3}d\overset{\cdot }{x}%
_{2}\wedge dx_{3}\wedge dt \\
&+&L_{\overset{\cdot }{x}_{1}}^{3}d\overset{\cdot }{x}_{1}\wedge
dx_{3}\wedge dt+L_{\overset{\cdot }{x}_{3}}^{1}d\overset{\cdot }{x}%
_{3}\wedge dx_{1}\wedge dt \\
&+&L_{\overset{\cdot }{x}_{2}}^{1}d\overset{\cdot }{x}_{2}\wedge
dx_{1}\wedge dt+L_{\overset{\cdot }{x}_{1}}^{2}d\overset{\cdot }{x}%
_{1}\wedge dx_{2}\wedge dt \\
&+& L_{\overset{\cdot }{x}_{1}}^{1}d\overset{\cdot }{x}_{1}\wedge
dx_{1}\wedge dt +L_{\overset{\cdot }{x}_{2}}^{2}d\overset{%
\cdot }{x}_{2}\wedge dx_{2}\wedge dt +L_{\overset{\cdot }{x}%
_{3}}^{3}d\overset{\cdot }{x}_{3}\wedge dx_{3}\wedge dt.
\end{eqnarray*}
Limiting it on vector fields $v=(\partial _{\Theta _{k}},\partial
_{d\Theta _{k}})$,
\begin{equation}
\label{eq5}\left( \text{rot}\,L\right) ^{k}=-\overset{\cdot }{
\lambda }^k,
\end{equation}
\[
L_{\overset{\cdot }{x}_{2}}^{3}-L_{\overset{ \cdot
}{x}_{3}}^{2}=2\lambda ^1, \qquad L_{\overset{\cdot
}{x}_{3}}^{1}-L_{\overset{ \cdot }{x}_{1}}^{3}=2\lambda ^2, \qquad
L_{\overset{\cdot }{x}_{1}}^{2}-L_{\overset{ \cdot
}{x}_{2}}^{1}=2\lambda ^3
\]
we get the Euler-Lagrange equation (\ref{eq4}).

\bigskip
\textbf{4.} Now we consider construction of the vector hamiltonian
$h^i$ for $L \in J^3(\pi)$. Rewrite the Griffiths 2-form $\psi$ as
\begin{eqnarray*}
\psi &=&L^i\omega _i+\lambda^i\Theta _i \\
&=&\left( L^{1}-\left( \lambda ^3\overset{\cdot }{x}_{2}-\lambda ^2%
\overset{\cdot }{x}_{3}\right) \right) dx_{1}\wedge dt \\
&+&\left( L^{2}-\left( \lambda ^1\overset{\cdot }{x}_{3}-\lambda ^3%
\overset{\cdot }{x}_{1}\right) \right) dx_{2}\wedge dt \\
&+&\left( L^{3}-\left( \lambda ^2\overset{\cdot }{x}_{1}-\lambda ^1%
\overset{\cdot }{x}_{2}\right) \right) dx_{3}\wedge dt+\lambda^idS_i \\
&=&-h^idx_i\wedge dt+\lambda ^idS_i.
\end{eqnarray*}
Here $dS_i=\varepsilon_{ijk}dx_j\wedge dx_k$ be a Pl\"ucker
coordinats of area element $dS$ spanned by vectors $dx_i$.

\textbf{Definition 1.} The vector field $\textbf{f}$ is called
conservative if
\[
\text{div}\;\textbf{f}=0.
\]

In other words, conservative vector field is divergence-free.

\textbf{Definition 2.} Phase trajectory $\textbf{x}(t)$ is called
Lagrange-stable if for all $t>0$ remains in some bounded domain of
phase space. Geometrically it means, that a phase flow (\ref{eq2})
should be divergence-free.

\textbf{Theorem 2.} The Lagrange-stable phase flow is hamiltonians.

\textbf{Proof.} We first calculate the exterior derivatives of
closed 2-forms $\psi$:
\[
d\psi =\left( \text{rot}\; \textbf{h}\right)^k \Theta _k\wedge
dt+\overset{\cdot }{ \lambda ^k}dt\wedge dS_{k}+
\text{div}\;\lambda^k \cdot dV=0.
\]
Then from $\text{div}\;\lambda ^k=0$ it follows that
\[
\overset{\cdot }{\lambda }=\text{rot}\;\textbf{h}.
\]
I.e. from hamiltonians point of view the set (\ref{eq2}) described
of dynamics of a generaliszed moments $\lambda$, which were defined
in (\ref{eq5}).

\bigskip
\textbf{5.} The base of deformation quantization of dynamical system
in  $J^2(\pi)$ is the Liouville theorem about preserved of the
volume $\Omega = dx_0 \wedge dx_1$ by phase flows. Geometrically it
means, that Lie derivative of the 2-form $\Omega $ along vector
field $X_H^1$ is zero: $\mathcal{L}_X \Omega = 0$. In other words if
$\{ g_t\}$ denotes the one parameter group symplectic
diffeomorphisms generated by vector fields $X_H^1$, then $g_t^\ast
\Omega = \Omega$ and the phase flow $\{g_t\}$ preserved the volume
form $\Omega$.

For extended this construction on $J^3(\pi)$ we consider the 3-form
of the phase space volume
\[
\Omega=dx_0 \wedge dx_1 \wedge dx_2.
\]

\textbf{Theorem 3.} The volume 3-form $ \Omega \in \Lambda ^3$
supposes existence two polyvector hamiltonians fields $X_H^1\in
\Lambda ^1$ and $X_H^2\in \Lambda ^2$.

\textbf{Proof.} By definition, put
\[
\mathcal{L}_{X} \Omega = X \rfloor d\Omega + d\left( {X\rfloor
\Omega} \right) = 0.
\]
Since $ \Omega \in \Lambda ^ 3 $, we see that $d\Omega = 0 $ and
\[
d\left( {X\rfloor \Omega}  \right) = 0.
\]
From Poincare's lemma it follows that form $X\rfloor \Omega $ is
exact, and
\[
X\rfloor \Omega = \Theta = d\textbf{H}.
\]
1) If $\textbf {X}_H^1 \in \Lambda ^1$, then $ \Theta \in \Lambda ^
2 $, $\textbf{H}=(\textbf{h} \cdot d\textbf{x}) \in \Lambda ^ 1 $.
Hamiltonian vector fields has the form
\begin{eqnarray}
X_H^1&=&(\text{rot}\,\textbf{h} \cdot \frac{\partial}{\partial
\textbf{x}})\label{eq6}\\
&=& \left(\frac{\partial h_2}{\partial x_1}- \frac{\partial
h_1}{\partial x_2} \right)\frac{\partial }{\partial x_0} + \left(
\frac{\partial h_0}{\partial x_2} - \frac{\partial h_2}{\partial
x_0} \right)\frac{\partial}{\partial x_1} + \left( \frac{\partial
h_1}{\partial x_0}-\frac{\partial h_0}{\partial x_1}
\right)\frac{\partial }{\partial x_2}.\notag
\end{eqnarray}

\noindent 2) If $X_H^2 \in \Lambda ^2$, then $ \Theta \in \Lambda ^1
$, $H \in \Lambda ^0 $ and we see already  hamiltonian bivector
fields
\begin{equation}
\label{eq7} X_{H}^{2} = \frac{{1}}{{2}}\left( {\frac{{\partial
H}}{{\partial x_{0}} } \cdot \frac{{\partial} }{{\partial x_{1}} }
\wedge \frac{{\partial }}{{\partial x_{2}} } + \frac{{\partial
H}}{{\partial x_{1}} } \cdot \frac{{\partial} }{{\partial x_{2}} }
\wedge \frac{{\partial} }{{\partial x_{0}} } + \frac{{\partial
H}}{{\partial x_{2}} } \cdot \frac{{\partial }}{{\partial x_{0}} }
\wedge \frac{{\partial} }{{\partial x_{1}} }} \right).
\end{equation}

More generalized (but scalar) construction was considered in
\cite{Dumachev}.

Poisson's bracket for vectorial hamiltonian (\ref{eq6}) has the form
\begin{eqnarray*}
\{ \textbf{h},G\}&=&X_H^1\rfloor dG\\
&=& \left(\frac{\partial h_2}{\partial x_1}- \frac{\partial
h_1}{\partial x_2} \right)\frac{\partial G}{\partial x_0} + \left(
\frac{\partial h_0}{\partial x_2} - \frac{\partial h_2}{\partial
x_0} \right)\frac{\partial G}{\partial x_1} + \left( \frac{\partial
h_1}{\partial x_0}-\frac{\partial h_0}{\partial x_1}
\right)\frac{\partial G}{\partial x_2},
\end{eqnarray*}
and dynamic equations is (\ref{eq2})
\begin{eqnarray}\label{eq8}
\overset{\cdot}{\textbf{x}}=\{ \textbf{h},\textbf{x} \}.
\end{eqnarray}

Poisson's bracket for bivector fields requires introduction two
hamiltonians
\begin{eqnarray}\label{eq9}
&&X_{H}^{2} \rfloor \left( {dF \wedge dG} \right) = \left\{ H,F,G
\right\}=\frac{1}{2}\left[  {\frac{{\partial H}}{{\partial x_{0}} }
\cdot \left( {\frac{{\partial F}}{{\partial x_{1}} }\frac{{\partial
G}}{{\partial x_{2} }} - \frac{{\partial F}}{{\partial x_{2}}
}\frac{{\partial G}}{{\partial x_{1}} }} \right)}\right. \notag \\
\\
&+&\left.  {\frac{{\partial H}}{{\partial x_{1}} } \cdot \left(
{\frac{{\partial F}}{{\partial x_{2}} }\frac{{\partial G}}{{\partial
x_{0} }} - \frac{{\partial F}}{{\partial x_{0}} }\frac{{\partial
G}}{{\partial x_{2}} }} \right) + \frac{{\partial H}}{{\partial
x_{2}} } \cdot \left( {\frac{{\partial F}}{{\partial x_{0}}
}\frac{{\partial G}}{{\partial x_{1} }} - \frac{{\partial
F}}{{\partial x_{1}} }\frac{{\partial G}}{{\partial x_{0}} }}
\right)}\right], \notag
\end{eqnarray}
such that dynamic equations (\ref{eq2}) has the Nambu form
\cite{Nambu}
\[
\overset{\cdot}{\textbf{x}}=\{ F,G,\textbf{x} \}.
\]

\bigskip
\textbf{Example.} Consider the dynamics of Frenet frame  with
constant curvature and torsion
\begin{eqnarray}\label{eq10}\left\{
\begin{array}{l}
\overset{\cdot}{x}=y\\
\overset{\cdot}{y}=z-x\\
\overset{\cdot}{z}=-y
\end{array}\right.
\end{eqnarray}
Lax representations for this set has the form
\[
\overset{\cdot}{A}=[A,B], \quad A= \left(
\begin{array}{ccc}
x&y&x\\
y&2z&y\\
x&y&x
\end{array}\right), \quad B= \left(
\begin{array}{ccc}
0&1&0\\
-1&0&-1\\
0&1&0
\end{array}\right)
\]
and gives us following invariants
\[
I_k=\frac{1}{k}\text{Tr}\textbf{A}^k,
\]
\[ I_1=x+z, \quad
I_2=\frac{1}{2}(x^2+y^2+z^2), \quad
I_3=\frac{1}{3}\left(x^3+\frac{3}{2}y^2(x+z)+z^3\right) ...
\]
Let $H_1=x+z$ and $H_2=\frac{1}{2}(2xz-y^2)$ - are the hamiltonians
of Frenet set, then
\[
I_1=H_1, \quad I_2=\frac{1}{2}H_1^2-H_2, \quad
I_3=\frac{1}{3}H_1\left(H_1^2-3H_2\right) ...
\]
The system (\ref{eq2}) is equivalent to system
\[
\overset{\cdot}{\textbf{x}}=\{ H_1,H_2,\textbf{x} \}
\]
with a Poisson bracket (\ref{eq9}).

For a finding of vectorial hamiltonian we write the differential
$\Psi=d\psi$ of Lagrange's 1-form for Frenet set (\ref{eq10}):
\[
\Psi=ydy\wedge dz+(z-x)dz\wedge dx-ydx\wedge dy,
\]
and, using homotopy formula, we get an expression for the vectorial
hamiltonians $h^i$ and vectorial lagranfians $L^i$:
\begin{eqnarray*}
\textbf{h}=\frac{1}{3}\left(
\begin{array}{l}
y^2+z^2-xz\\
-y(x+z)\\
y^2+x^2-xz
\end{array}\right), \qquad
\textbf{L}=\left(
\begin{array}{l}
z\overset{\cdot}{y}-y\overset{\cdot}{z}-h_1\\
x\overset{\cdot}{z}-z\overset{\cdot}{x}-h_2\\
y\overset{\cdot}{x}-x\overset{\cdot}{y}-h_3\\
\end{array}\right).
\end{eqnarray*}

\end{document}